\documentclass[12pt]{article}

\usepackage{amsmath,amsthm,amssymb}
\usepackage{xypic}

\newcommand{\namedthm}[2]{\theoremstyle{plain}
   \newtheorem*{thm#1}{#1}\begin{thm#1}#2\end{thm#1}}

\newtheorem{lemma}{Lemma}[section]
\newtheorem{theorem}{Theorem}
\newtheorem{corollary}{Corollary}
\newtheorem*{definition*}{Definition}
\newtheorem{proposition}{Proposition}

\def\R{\mathbb{R}}
\def\C{\mathbb{C}}
\def\Q{\mathbb{Q}}
\def\Z{\mathbb{Z}}
\def\O{\mathcal{O}}
\def\T{\mathbb{T}}
\def\U{\mathcal{U}}
\def\eps{\varepsilon}
\def\qed{$\blacksquare$}

\DeclareMathOperator{\Orb}{Orb}
\DeclareMathOperator{\rank}{rank}
\DeclareMathOperator{\length}{length}
\DeclareMathOperator{\Tr}{Tr}
\DeclareMathOperator{\Log}{Log}

\title{On the $S$-Euclidean minimum of an ideal class}
\author{Kevin J. McGown}

\date{}

\begin{document}

\maketitle

\vspace{-3ex}

\begin{abstract}
We show that the $S$-Euclidean minimum of an ideal class is a rational number, generalizing a result of Cerri.
We also give some corollaries which explain the relationship of our results with
Lenstra's notion of a norm-Euclidean ideal class and the conjecture of Barnes and Swinnerton-Dyer on quadratic forms.
The proof is self-contained but uses ideas from ergodic theory and topological dynamics,
particularly those of Berend.
\end{abstract}


\section{Introduction}\label{S:intro}

Let $K$ be a number field of degree $n=r_1+2r_2$.  Let $S$ be a finite set of primes containing the infinite primes $S_\infty$.
Let $\O_S$ denote the $S$-integers of $K$, let $\U_S$ denote the $S$-units of $K$,
and let $N_S$ denote the $S$-norm map.\footnote{In cases where ``$S$'' is dropped from the notation
this will mean we are using $S=S_\infty$.}
Recall that the $S$-norm of a number $\xi\in K$ is defined as
$N_S(\xi)=\prod_{v\in S}|\xi|_v$
and the $S$-norm of an ideal $\mathfrak{a}\subseteq\O_S$ is defined as $N_S(\mathfrak{a})=|\O_S/\mathfrak{a}|$.

For an ideal $\mathfrak{a}\subseteq\O_S$ and an element $\xi\in K$, we define
$$
  m_{\mathfrak{a}}^S(\xi)
  =
  \frac{1}{N_S(\mathfrak{a})}
  \;
  \inf_{\gamma\in\mathfrak{a}}
  N_S(\xi-\gamma)
  \,,
$$
and we set
$$
  M_{\mathfrak{a}}^S
  =
  \sup_{\xi\in K}
  m_{\mathfrak{a}}^S(\xi)
  \,.
$$
Notice that $m_{\mathfrak{a}}^S(\xi)$ depends upon the ideal $\mathfrak{a}$,
but that $M_{\mathfrak{a}}^S$ only depends upon the ideal class $[\mathfrak{a}]$;
this follows easily from the fact $\mathfrak{a}=\gamma\mathfrak{b}$ implies $m_\mathfrak{a}(\xi)=m_\mathfrak{b}(\xi\gamma^{-1})$
for nonzero $\gamma\in K$.
We call $M_\mathfrak{a}^S$ the $S$-Euclidean minimum of the ideal class~$[\mathfrak{a}]$.

One easily verifies that $m_\mathfrak{a}^S(\xi)\in\Q$ for all $\xi\in K$
since the infimum is being taken over a discrete subset of $\Q$.\footnote{
Given $\xi\in K$ there exists $d\in\Z^+$ such that $d\xi\in\O$
and hence $\{N_S(\xi-\gamma)\}_{\gamma\in\mathfrak{a}}$ is contained in $ N_S(d)^{-1}\Z$.
}
However, it is not by any means clear whether $M_\mathfrak{a}^S$ is rational or not.
When $K=\Q(\sqrt{d})$, $d>0$, and $S=S_\infty$, the statement
$M_\mathfrak{a}^S\in \Q$ is equivalent to a classical conjecture of Barnes and Swinnerton-Dyer, which is still unresolved.
Our aim is to prove the following:

\begin{theorem}\label{T:1}
If $\#S\geq 3$, then $M_{\mathfrak{a}}^S\in \Q$.
\end{theorem}

Cerri proved (see~\cite{cerri}) that
$M_\mathfrak{a}^S\in\Q$ when $S=S_\infty$, $\# S\geq 3$, and $\mathfrak{a}=\O_K$;
our theorem generalizes his result in two directions --- to the $S$-integral setting
and to ideal classes.

As we will discuss in \S\ref{S:euclidean}, the quantity $M_\mathfrak{a}^S$
is important in the study of norm-Euclidean ideal classes.
\emph{A~priori}, it is possible that there are norm-Euclidean ideal classes with
\mbox{$M_\mathfrak{a}^S=1$}.
However, our results lead to:
\begin{corollary}\label{C:1}
If $\#S\geq 3$, then an ideal class $[\mathfrak{a}]$ of $\O_S$ is norm-Euclidean if and only if
$M_\mathfrak{a}^S<1$.
\end{corollary}

In \S\ref{S:forms} we discuss the relationship of the quantity $M_\mathfrak{a}^S$
with the conjecture of Barnes
and Swinnerton-Dyer.
In the case of $K=\Q(\sqrt{d})$, $d>0$ one has $\#S_\infty=2$ and hence our result
gives the following:
\begin{corollary}\label{C:2}
``The conjecture of Barnes and Swinnerton--Dyer holds for fundamental discriminants
if we invert a single prime.''
\end{corollary}

Strictly speaking, the previous two corollaries will follow from a slight refinement of Theorem~\ref{T:1}
which we describe after introducing the requisite notation (see Theorem~\ref{T:2}).


\section{Main idea and setup}\label{S:main.idea}
From the definition one sees that $m_{\mathfrak{a}}^S$ can be viewed as a function $K\to \R_{\geq 0}$
as well as a function $K/\mathfrak{a}\to \R_{\geq 0}$.
As there should be no confusion, we will
denote both of these functions by
$m_{\mathfrak{a}}^S$.

The $S$-units $\U_S$ act on $K/\mathfrak{a}$ by multiplication
and the function $m_{\mathfrak{a}}^S:K/\mathfrak{a}\to \R_{\geq 0}$
is invariant under this action.  The main idea is to embed $K/\mathfrak{a}$
into a compact metric group $\mathbb{T}$ where $\U_S$ still acts and $m_{\mathfrak{a}}^S$
extends naturally to an upper semi-continuous function on $\mathbb{T}$.
In this setting we will be able to study the action of the units
from the point of view of ergodic theory and topological dynamics.

We embed $K$ diagonally into the product of its completions at the primes in $S$.
We will write $K\subseteq\prod_{v\in S} K_v=:\overline{K}_S$.
The function $N_S:K\to\R_{\geq0}$ extends to a continuous function $N_S:\overline{K}_S\to\R_{\geq 0}$,
and this allows us to define $m_{\mathfrak{a}}^S(\xi)$ for any $\xi\in\overline{K}_S$.
It follows that $m_{\mathfrak{a}}^S:K\to\R_{\geq 0}$ extends to 
an upper semi-continuous function $m_{\mathfrak{a}}^S:\overline{K}_S\to\R_{\geq 0}$.
Finally, we define
$$
  \overline{M}_\mathfrak{a}^S:=\sup_{\xi\in\overline{K}_S}m_\mathfrak{a}^S(\xi)
  \,.
$$
We call $\overline{M}_\mathfrak{a}^S$ the $S$-inhomogeneous minimum of the ideal class
$[\mathfrak{a}]$.

The embedding $K\subseteq \overline{K}_S$ induces an embedding
$K/\mathfrak{a}\subseteq\overline{K}_S/\mathfrak{a}=:\T$,
and $m_{\mathfrak{a}}^S$ induces an upper semi-continuous function
$m_{\mathfrak{a}}^S:\T\to\R_{\geq 0}$.  Since $\T$ is compact,
this tells us that there exists $\xi_0\in \overline{K}_S$ such that
$m_\mathfrak{a}^S(\xi_0)=\overline{M}_\mathfrak{a}^S$.
If we can show that there exists an element $\xi_0\in K$ such that 
$m_\mathfrak{a}^S(\xi_0)=\overline{M}_\mathfrak{a}^S$, then it would follow that
$M_\mathfrak{a}^S=\overline{M}_\mathfrak{a}^S\in\Q$.
Indeed, we will prove the following result from which
Theorem~\ref{T:1} follows.
\begin{theorem}\label{T:2}
Suppose $\#S\geq 3$.
Then there exists an element $\xi_0\in K$ such that 
$m_\mathfrak{a}^S(\xi_0)=\overline{M}_\mathfrak{a}^S$.
\end{theorem}

Although many aspects of the proof of Theorem~\ref{T:2}
are motivated by ideas in
topological dynamics and ergodic theory
(particularly those of Berend),
our account will be largely
self-contained.  In fact, except for a couple small lemmas,
the only outside results we appeal to
are standard theorems in number theory.
However, we should mention that
much has been gleaned from studying the
papers~\cite{furstenberg,berend:multi,berend:minimal,berend:compact,cerri}.

Before proceeding to the proof of Theorem~\ref{T:2}, we will
discuss some applications.


\section{Euclidean Ideal Classes}\label{S:euclidean}
Lenstra introduced the following definition:
We call an ideal class $[\mathfrak{a}]$ of $\O_S$
norm-Euclidean if for every $\xi\in K$ there exists $\gamma\in\mathfrak{a}$ such that
$N_S(\xi-\gamma)<N_S(\mathfrak{a})$.
(Recall that our norms are defined to be positive.)
Notice that that if we take $\mathfrak{a}=(1)$, then this reduces to the usual definition
of the ring $\O_S$ being norm-Euclidean.
One important fact is that if an ideal class $[\mathfrak{a}]$ is norm-Euclidean, then
it generates the class group of $\O_S$; in particular, the existence of a norm-Euclidean
ideal class implies that the class group is cyclic.  (See~\cite{lenstra} for more details.)

It is clear that $M_\mathfrak{a}^S<1$ implies that $[\mathfrak{a}]$ is
norm-Euclidean, and that $M_\mathfrak{a}^S>1$ implies that $[\mathfrak{a}]$ is
not norm-Euclidean.  In the case $M_\mathfrak{a}^S=1$,
one cannot draw any immediate conclusion.
However, in light of Theorem~\ref{T:2}, provided $\#S\geq 3$, the condition
$M_\mathfrak{a}^S=1$ always implies that $[\mathfrak{a}]$ is not norm-Euclidean;
indeed, in this case,
Theorem~\ref{T:2} implies
that there exists $\xi_0\in K$
such that $N_S(\xi_0-\gamma)\geq N_S(\mathfrak{a})$ for all $\gamma\in\mathfrak{a}$.
This establishes Corollary~\ref{C:1}.

Define the open neighborhoods
$V_t:=\{\xi\in\overline{K}_S\mid N_S(\xi)<t\}$.
Lenstra points out that $[\mathfrak{a}]$ is norm-Euclidean
if and only if $K\subseteq \mathfrak{a}+V_{N_S(\mathfrak{a})}$.
We quote~\cite{lenstra} (using our notation):
``It seems that in all cases in which this condition is known to be satisfied we actually have
$\overline{K}_S=\mathfrak{a}+V_{N_S(\mathfrak{a})}$.
It is unknown whether both properties are in fact
equivalent.''\footnote{He then goes on to state the only known result
in this direction.
It is not important to us here as it pertains to the case where $\#S\leq 2$.}
We completely answer this question when $\#S\geq 3$ (in the number field case) with the following:
\begin{corollary}\label{C:euclid1}
Suppose $\#S\geq 3$.  Then we have
$K\subseteq \mathfrak{a}+V_{N_S(\mathfrak{a})}$ if and only if
$\overline{K}_S=\mathfrak{a}+V_{N_S(\mathfrak{a})}$.
\end{corollary}

\noindent\textbf{Proof.}
One direction of the result is obvious.
To prove the other direction, suppose 
$K\subseteq \mathfrak{a}+V_{N_S(\mathfrak{a})}$;
in other words, $[\mathfrak{a}]$ is norm-Euclidean.
In light of Theorem~\ref{T:2} and Corollary~\ref{C:1} we see that
$\overline{M}_\mathfrak{a}^S=M_\mathfrak{a}^S<1$.
It follows that for every $\xi\in\overline{K}_S$ there exists $\gamma\in\mathfrak{a}$
such that $N_S(\xi-\gamma)<N_S(\mathfrak{a})$;
this proves
$\overline{K}_S\subseteq\mathfrak{a}+V_{N_S(\mathfrak{a})}$.
\qed

\vspace{1ex}

In light of the discussion in~\cite{lenstra:intelligencer2},
we now immediately obtain the following additional result:
\begin{corollary}\label{C:euclid2}
The question of whether $[\mathfrak{a}]$ is norm-Euclidean
is decidable when $\#S\geq 3$.
\end{corollary}
Readers interested in reading more regarding the Euclidean algorithm in number fields
should consult the excellent expository article~\cite{lemmermeyer:survey}.


\section{The Conjecture of Barnes and Swinnerton-Dyer}\label{S:forms}

Let $f(x,y)=ax^2+bxy+cy^2$ with $a,b,c\in\Z$ be a binary quadratic form
with discriminant $\Delta=b^2-4ac>0$.
For ease of exposition, we will henceforth write \emph{form} to mean \emph{binary quadratic form}.
For a form $f$ and a point $P\in\Q^2$, we define
$$
  m_{f}(P)
  =
  \inf_{Q\in\Z^2}
  |f(P-Q)|
  \,,
$$
and we set
$$
  M_{f}
  =
  \sup_{P\in\Q^2}
  m_{f}(P)
  \,,
  \quad
  \overline{M}_{f}
  =
  \sup_{P\in\R^2}
  m_{f}(P)  
  \,.
$$
Since $M_{\lambda f}=|\lambda|\, M_f$ and
$\overline{M}_{\lambda f}=|\lambda|\, \overline{M}_f$ for all $\lambda\in\Z$,  
we will only consider forms where $\gcd(a,b,c)=1$, which are known as
primitive forms.
Barnes and Swinnerton-Dyer conjecture
(see~\cite{barnes})
that there exists a point $P_0\in\Q^2$ such that
$m_f(P_0)=M_f=\overline{M}_f$; in particular,
$M_f\in \Q$.\footnote{They also conjecture that the minimum is so-called attained and isolated,
but we will ignore this part of the conjecture for the purposes of
this investigation.}

Fix a fundamental discriminant $\Delta>0$.
Let $K=\Q(\sqrt{\Delta})$ be the real quadratic field of discriminant $\Delta$
having ring of integers $\O$.
Let $\mathfrak{a}$ be an ideal of $\O$ with $\Z$-basis $\{\alpha_1,\alpha_2\}$.
We can associate to $\mathfrak{a}$ the form of discriminant $\Delta$ given by
$$
  \frac{1}{N(\mathfrak{a})}(\alpha_1 x+\alpha_2 y)(\overline{\alpha}_1 x+\overline{\alpha}_2 y).
$$
In fact, every primitive form of discriminant $\Delta$
arises in this way.\footnote{One can extend the correspondence to include forms with non-fundamental discriminants
by considering orders other than the full ring of integers, but
in this paper we are content to restrict ourselves to forms with fundamental discriminants.}
See~\cite{buell} for a classical treatment of this correspondence or~\cite{bhargava} for a more modern treatment.

The conjecture of Barnes and Swinnerton-Dyer (as stated above)
for fundamental discriminants is equivalent to the statement:
Given an ideal class $[\mathfrak{a}]$ in a real quadratic field
$K$, there exists $\xi_0\in K$ such that $m_\mathfrak{a}^{S_\infty}(\xi_0)=\overline{M}_\mathfrak{a}^{S_\infty}$.
Although we cannot prove this statement, since $\#S_\infty=2$ in the case where
 $K=\Q(\sqrt{d})$, $d>0$,
we can prove the analogous statement when $S=S_\infty\cup\{\mathfrak{p}\}$
where $\mathfrak{p}$ is any (finite) prime of $K$.  This follows from Theorem~\ref{T:2}
and is the content of Corollary~\ref{C:2}.


\section{Preliminary Results}\label{S:prelim}

In this section we give a brief justification for the facts claimed in \S\ref{S:main.idea}
and derive a couple other basic results.
The hurried reader who is willing to
refer back to this section as necessary
may skip to \S\ref{S:next}.

Observe that $\overline{K}_S$ is a locally compact abelian group.
It is also a complete metric space
with metric $d(\alpha,\beta)=\max_{v\in S} |\alpha_v-\beta_v|_v$. 
The fact that $N_S(\xi)=\prod_{v\in S}|\xi_v|_v$ is continuous on $\overline{K}_S$
follows immediately from the fact that each
$|\cdot|_v:K_v\to\R$ is continuous.

To show that  $\O_S$ is discrete in $\overline{K}_S$,
it suffices to show that $\{0\}$ is open in the subspace topology on $\O_S$.
The set $V=\{\alpha\in\overline{K}_S\mid N_S(\alpha)<1\}$
is open in $\overline{K}_S$ since $N_S$ is continuous, and, moreover,
$V\cap \O_S=\{0\}$.
Since $\O_S$ is discrete in $\overline{K}_S$, so is $\mathfrak{a}$.
It now follows from generalities that
$\T$ is a locally compact Hausdorff space.
In fact, one can show that the metric on $\overline{K}_S$ induces a metric on $\T$
in the usual manner.

The only fact that remains to be justified is that $\T$ is compact.
For this, we will need the following standard result from algebraic number theory
(see, for example, \cite{cassels}).


%
%

\namedthm{Strong Approximation Theorem}
{
Suppose we are given a finite set of primes $T$, elements $\alpha_v\in K_v$ for each $v\in T$, and
a prime $w\notin T$.
Then for each $\eps>0$, there exists a number $\beta\in K$ such that $|\alpha_v-\beta|_v<\eps$
for all $v\in T$ and $|\beta|_v\leq 1$ for all $v\notin T$ with $v\neq w$.
}

We mention in passing that applying the previous result with $T=S$
tells us that $K$ is dense in $\overline{K}_S$,
which explains the notation.
In what follows we write $S_0$ for the finite primes in $S$ so that $S=S_\infty\cup S_0$.

\begin{lemma}\label{L:approx.help}
Let $(\alpha_v)_{v\in S}\in \overline{K}_S$.
Then there exists $\gamma\in\mathfrak{a}$ such that \mbox{$v(\alpha_v-\gamma)\geq 0$}
for all $v\in S_{0}$.
\end{lemma}

\noindent\textbf{Proof.}
By the Strong Approximation Theorem, there exists $\gamma\in K$ such that
$v(\alpha_v-\gamma)\geq 0$ for all 
$v\in S_{0}$,
$v(\gamma)\geq v(\mathfrak{a}\cap\O)$ for all finite $v$ dividing $\mathfrak{a}\cap\O$,
and $v(\gamma)\geq 0$ for all other $v$.
This is possible since $\mathfrak{a}\cap\O$ is not divisible by any primes in $S$.
One checks that this choice of $\gamma$ works.
\qed

\begin{lemma}
Let $\mathcal{F}$ be a fundamental domain for $\overline{K}_{S_\infty}/(\mathfrak{a}\cap\O)$.
(Note that $\overline{K}_{S_\infty}\simeq\R^{r_1}\times\C^{r_2}$ is the usual
Minkowski space and $\mathfrak{a}\cap\O$ is an ideal of $\O$.)
Each element of $\T$ has a unique representative in
$\mathcal{F}\times\prod_{v\in S_0}\O_v$.
\end{lemma}

\noindent\textbf{Proof.}
Let $(\alpha_v)\in\overline{K}_S$.
Using Lemma~\ref{L:approx.help},
choose $\gamma\in\mathfrak{a}$ such that $v(\alpha_v-\gamma)\geq 0$ for all $v\in S_0$.
Choose $a\in\O\cap\mathfrak{a}$ such that $(\alpha_v-\gamma-a)_{v\in S_\infty}\in\mathcal{F}$.
Then $\alpha_v-\gamma-a\in\O_v$ for all $v\in S_0$.
Set $\beta_v=\alpha_v-\gamma-a\in\mathcal{F}\times\prod_{v\in S_0}\O_v$.
Then $(\alpha_v)-(\beta_v)=\gamma+a\in\mathfrak{a}$.

Now we show uniqueness.
Suppose $(\alpha_v)=(\beta_v)+\delta$ for some $\delta\in\mathfrak{a}$
with $(\alpha_v),(\beta_v)\in\mathcal{F}\times\prod_{v\in S_0}\O_v$.
Then
$v(\delta)\geq 0$ for all $v\in S_0$
which implies
$\delta\in\O\cap\mathfrak{a}$.
Since
$(\alpha_v)_{v\in S_\infty},(\beta_v)_{v\in S_\infty}\in\mathcal{F}$ and $\delta\in\O\cap\mathfrak{a}$,
we have $(\alpha_v)_{v\in S_\infty}=(\beta_v)_{v\in S_\infty}$ which implies $\delta=0$.
\qed

\begin{lemma}\label{L:compact}
$\T$ is compact.
\end{lemma}

\noindent\textbf{Proof.}
By the previous lemma, $\T$ is the image under the natural projection of the compact set
$\overline{\mathcal{F}}\times\prod_{v\in S_{0}}\O_v$.
\qed

\vspace{1ex}
We conclude this section with another simple result
that is a consequence of Strong Approximation
which will prove useful in the sequel.
For each $w\in S$, we can view $K_w$ as a subset of $\overline{K}_S$ by sending
the element $x\in K_w$ to the vector $\xi\in\overline{K}_S$ where $\xi_w=x$ and $\xi_v=0$ for $v\neq w$;
that is, the image of $K_w$ in $\overline{K}_S$ is zero outside the $w$-component.

\begin{lemma}\label{L:prelim.help}
For each $w\in S$, we have $K_w+\mathfrak{a}$ is dense in $\overline{K}_S$.
In particular, there are no
proper closed subgroups of $\overline{K}_S$ containing both $K_w$ and $\mathfrak{a}$.
\footnote{Keep in mind that $K_w$ is embedded into one component and
$\mathfrak{a}$ is embedded diagonally.}
\end{lemma}

\noindent\textbf{Proof.}
Let $\xi=(\xi_v)_{v\in S}\in\overline{K}_S$.
Fix $\eps>0$.  
Using Strong Approximation, 
choose $\gamma\in K$ such that
$|\gamma-\xi_v|_v<\eps$ for all $v\in S$ with $v\neq w$,
$|\gamma|_v<\eps$ for all finite $v$ dividing $\mathfrak{a}\cap\O$,
and $|\gamma|_v\leq 1$ for all $v\notin S$.
When $\eps>0$ is small enough,
this implies $\gamma\in\mathfrak{a}$.
Additionally,
$\xi-\gamma=(\xi_v-\gamma)_{v\in S}$ is $\eps$-close to
$\beta:=\xi_w-\gamma\in K_w\subseteq\overline{K}_S$.
It follows that $\xi$ is $\eps$-close to $\beta+\gamma\in K_w+\mathfrak{a}$.
Since $\eps>0$ was arbitrary, there are points in $K_w+\mathfrak{a}$ arbitrarily close to $\xi$.
\qed


\section{Outline of the Proof}\label{S:next}

Given an element $\xi\in\overline{K}_S$, we will write $[\xi]$ for the class of $\xi$ in $\T$;
namely, $[\xi]=\pi(\xi)$ where
$\pi:\overline{K}_S\to\T$ is the natural projection map.
We begin with two lemmas concerning the orbit structure of the action of
$\U_S$ on $\T$.

\pagebreak 

\begin{lemma}\label{L:1}
For $\xi\in\overline{K}_S$, the following are equivalent:
\begin{enumerate}
\item
$\xi\in K$
\item
$[\xi]\in\mathbb{T}_{\text{tors}}$
\item
$\Orb([\xi])$ is finite
\item
$[\xi]$ is an isolated point of $\overline{\Orb([\xi])}$
\end{enumerate}
\end{lemma}

\noindent\textbf{Proof.}
First we show $(1)\Rightarrow(2)\Rightarrow(3)$.
If $\xi\in K$, there exists $n\in\Z^+$ such that $n\xi\in\mathfrak{a}$
and hence $[\xi]\in\T_{\text{tors}}$.
In this case $u\xi\in(1/n)\mathfrak{a}$ for all $u\in\U_S$ and therefore
$\Orb([\xi])\subseteq\pi((1/n)\mathfrak{a})$, which is a finite subgroup
of $K/\mathfrak{a}$.

Now we show $(3)\Rightarrow(1)$.  Suppose $u[\xi]=u'[\xi]$ in $\T$ with $u\neq u'\in\U_S$.
Then there exists $\alpha\in\mathfrak{a}$ such that $u\xi=u'\xi+\alpha$ in $\overline{K}_S$.
It follows that $u\xi_v=u'\xi_v+\alpha$ in $K_v$ for all $v\in S$.
We conclude that $\xi_v=\alpha/(u-u')\in K$ for all $v$ and therefore
$\xi\in K$.

Now we show $(3)\Leftrightarrow(4)$.
For convenience of notation let $A=\overline{\Orb([\xi])}$.
The set $A$ is a closed subset of $\T$ and therefore compact
(see Lemma~\ref{L:compact}).
It is now easy to see that for
$\xi\in\overline{K}_S$
one has:
$[\xi]$ is isolated in $A$ iff
$\Orb([\xi])$ is discrete in $A$ iff $\Orb([\xi])$ is finite.
\qed

\begin{lemma}
Let $\xi\in\overline{K}_S\setminus K$.
Then the map $\U_S\to\Orb([\xi])$ given by $u\mapsto u[\xi]$
is a bijection.
\end{lemma}

\noindent\textbf{Proof.}
This follows immediately from the proof of $(3)\Rightarrow(1)$ 
in the previous lemma.
\qed

\vspace{1ex}

The next lemma is easily deduced, but essential.
It constitutes the natural generalization
of an important observation of Cerri.
In fact, we employ the group $\T=\overline{K}_S/\mathfrak{a}$ precisely so that the following
result will go through in our setting:

\begin{lemma}\label{L:2}
The set 
$\{[\xi]\in\mathbb{T}\mid m_\mathfrak{a}^S(\xi)=\overline{M}_\mathfrak{a}^S\}$
is a nonempty closed $\mathcal{U}_S$-invariant subset of $\mathbb{T}$.

\end{lemma}

\noindent\textbf{Proof.}
This follows from the fact that $m_\mathfrak{a}^S$ is a $\U_S$-invariant, upper semi-continuous
function defined on the compact set $\mathbb{T}$.
\qed

\begin{theorem}\label{T:3}
Suppose $\#S\geq 3$.
Then every nonempty closed $\mathcal{U}_S$-invariant subset of $\mathbb{T}$ contains torsion elements.
\end{theorem}

If we can prove Theorem~\ref{T:3}, then Theorems~\ref{T:2} and~\ref{T:1} immediately follow
in light of Lemmas~\ref{L:1} and~\ref{L:2}.
The proof requires the next three propositions
whose proofs we postpone until Sections~\ref{S:Proof2}, \ref{S:Proof1}, and~\ref{S:CM} respectively.

\begin{definition*}
We refer to a nonempty $\U_S$-invariant closed subset of $\T$ which is minimal with respect
to set inclusion as a $\U_S$-minimal set.
\end{definition*}

Observe that by Zorn's Lemma, every nonempty $\U_S$-invariant closed subset
of $\T$ contains a $\U_S$-minimal set.

\begin{proposition}\label{P:3}
Let $M$ be $\U_S$-minimal subset of $\T$.
Then $M-M$ is a proper subset of $\T$.
\end{proposition}

Recall that a CM-field is a totally complex quadratic extension of a totally real field.
Let $K^+$ denote the maximal totally real subfield of $K$.
In the case where $K$ is a CM-field we have $[K:K^+]=2$.

\begin{proposition}\label{P:2}
Suppose $K$ is not a CM-field
or $S$ contains a finite prime that splits in $K/K^+$.
Let $N$ be a closed $\U_S$-invariant subset of $\T$
that contains $0$ as a non-isolated point.
If $\#S\geq 3$, then $N=\T$.
\end{proposition}

\begin{proposition}\label{P:4}
If Theorem~\ref{T:3} holds except in the case where $K$ is a CM-field and
distinct primes in $S$ lie over distinct primes in $K^+$,
then it holds in all cases.
\end{proposition}

\noindent\textbf{Proof of Theorem~\ref{T:3}.}
Let $M$ be a $\U_S$-minimal subset of $\T$.
Moreover, assume that $M$ contains no torsion elements.
Then $N=M-M$ is a closed $\U_S$-invariant subset of $\T$.

We show that $N$ contains $0$ as a non-isolated point.
Pick $[\xi]\in M$.  Then $M=\overline{\Orb([\xi])}$.
By Lemma~\ref{L:1}, $[\xi]$ must be non-isolated and therefore
there is a sequence $u_n\in\U_S$ with the $u_n$ distinct
and $u_n [\xi]\to [\xi]$.  Without loss of generality, we may assume that
$u_n\neq 1$ for all $n$.  Now observe that $u_n[\xi] - [\xi]$
is a sequence of nonzero points in $N$ converging to $0$.

By Proposition~\ref{P:4} it suffices to prove the theorem
in the situation where Proposition~\ref{P:2} applies.
We invoke Proposition~\ref{P:2}
and conclude that $N=\T$.
This contradicts Proposition~\ref{P:3}.
Thus $M$ contains torsion elements.
\qed


\section{Character Theory and Ergodicity}

Before turning to the proofs of Propositions~\ref{P:3},\ref{P:2},\ref{P:4}, we need a few lemmas
which are consequences of the study of the character theory of
$\overline{K}_S$ and $\T$.
We have tried to assume a minimal amount of background,
giving the appropriate definitions and
stating the necessary facts, but some familiarity with the
duality theory of locally compact abelian groups
(and local fields in particular)
will be helpful in this section.

Let $G$ be a locally compact abelian group.
A (unitary) character of $G$ is a continuous group homomorphism
$\chi:G\to S^1$.
(We will always view $S^1$ as the unit circle inside $\C$.)
The Pontryagin dual of $G$, denoted by $G^\vee$, is the
(abelian) multiplicative group consisting of all the characters of $G$.
It is locally compact when endowed with the
topology of uniform convergence on compact sets.

We are interested in the characters of $\T$.
However, since any character of $\T$ may be viewed as a character of $\overline{K}_S$
that is trivial on $\mathfrak{a}$, we will first consider characters of $\overline{K}_S$.
(Note that another way to view the previous sentence is that we have an injection
$\T^\vee\hookrightarrow\overline{K}_S^\vee$.)
The group $\overline{K}_S$ is self-dual since it is a product of local fields.
We now construct an explicit nontrivial character of $\overline{K}_S$
that will facilitate subsequent arguments.


\subsection{Constructing a character of $\overline{K}_S$}

For each $v\in S$, one can define a nontrivial local character $\phi_v:K_v\to S^1$
in a natural way.
If $v$ is real, we set $\phi_v(x)=e^{2\pi i x}$, and if $v$ is complex,
we set $\phi_v(z)=e^{2\pi i(z+\overline{z})}$.
In the case where $v$ is a finite prime, we define
$\phi_v(\alpha)$ to be the exponential of $2\pi i$ times the
the ``polar part'' of $\Tr_{K_v/\Q_p}(\alpha)$.\footnote{Here $p$ is the rational prime lying under $v$,
and the polar part of an element of $\Q_p$ is
the element of $\Q/\Z$ defined
by the (non-unique) decomposition
$\Q_p=\Z[1/p]+\Z_p$.}
Then $\psi=\prod_{v\in S}\phi_v$ is a non-trivial character of $\overline{K}_S$ and
we have the explicit isomorphism $\overline{K}_S\to\overline{K}^\vee_S$
given by $\xi\mapsto\psi_\xi$; here $\psi_\xi(\eta)=\psi(\xi\eta)$.

Since $\O_S$ is a proper closed subgroup of $\overline{K}_S$
there is a nonzero character $\phi$
of $\overline{K}_S$ that is trivial on $\O_S$.
By duality, $\phi=\psi_\rho$ for some $\rho\in\overline{K}_S$.
Therefore
\begin{equation}\label{E:phi}
  \phi(\xi)=\prod_{v\in S}\phi_v(\rho_v\xi_v)
  \,.
\end{equation}
As before we have an explicit isomorphism
$\overline{K}_S\to\overline{K}^\vee_S$ given by $\xi\mapsto\phi_\xi$ where $\phi_\xi(\eta)=\phi(\xi\eta)$.
To completely justify this, one should check that $\rho_v\neq 0$ for all $v\in S$.
Suppose it were the case that $\rho_w=0$ for some $w$.
Then we would have that $\ker\phi$ contains both $K_w$ and $\O_S$
which implies that $\phi$ is the trivial character
(by an application Lemma~\ref{L:prelim.help} with $\mathfrak{a}=\O_S$),
a clear contradiction.  We will not need to determine the $\rho_v$; it will be enough to know that they
are all nonzero.

\textbf{For the remainder of the paper, $\phi$ will refer to this particular fixed character of $\overline{K}_S$.
Likewise, the notations $\phi_\xi$ and $\phi_v$ will refer to the characters constructed here.}
Notice that $\phi$ depends upon the number field $K$ and set of primes $S$,
but it does not depend upon the choice of $\mathfrak{a}$.


\subsection{The dual of $\T$}

\begin{definition*}
Given a subset $E\subseteq\overline{K}_S$, we define
its complement:
$$
  E^\perp=\{\xi\in\overline{K}_S\mid \phi(\xi E)=1\}
  \,.
$$
\end{definition*}
It is an easy exercise to show
that if $E$ is a subgroup (or $\O_S$-submodule) of $\overline{K}_S$, then so is $E^\perp$.

\begin{lemma}\label{L:dual}
The map $\mathfrak{a}^\perp\to\T^\vee$ given by $\alpha\mapsto\phi_\alpha$
is an isomorphism of topological groups.
\end{lemma}

\noindent\textbf{Proof.}
Every character of $\T$ may be viewed as a character of $\overline{K}_S$ that is
trivial on $\mathfrak{a}$.  Every character of $\overline{K}_S$ is of the form
$\phi_\xi$ for some $\xi\in\overline{K}_S$.  Finally, a character $\phi_\xi$
is trivial on $\mathfrak{a}$
if and only if $\xi\in\mathfrak{a}^\perp$.
\qed

\vspace{1ex}

\noindent
To make the previous result useful, one would like a better description of $\mathfrak{a}^\perp$.

\begin{lemma}
If $\mathfrak{b}$ is a fractional ideal of $\O_S$, then so is $\mathfrak{b}^\perp$.
Moreover\\ $\mathfrak{b}^\perp=\mathfrak{b}^{-1}\O_S^\perp$.
\footnote{Here $\O_S^\perp$ plays the role of the inverse different $\mathfrak{D}^{-1}$.
In particular, when $S=S_\infty$ one can take $\phi=\prod_{v\in S}\phi_v$ and we have $\O_S^\perp=\mathfrak{D}^{-1}$.}
\end{lemma}

\noindent\textbf{Proof.}
Since $\phi(\O_S)=1$ we have $\O_S\subseteq\O_S^\perp$.  We show
that $\O_S^\perp/\O_S$ is finite.
First, since $\O_S^\perp$ is dual to the compact group $\overline{K}_S/\O_S$
(by the previous lemma)
we know that $\O_S^\perp$ is discrete.  It follows that $\O_S^\perp/\O_S$
is a discrete subspace of the compact space $\overline{K}_S/\O_S$ and therefore finite.
If we set $d=|\O_S^\perp/\O_S|$, then this gives $d\O_S^\perp\subseteq \O_S$
and therefore $\O_S^\perp$ is contained in $K$.
In light of previous comments, it now follows easily that
$\O_S^\perp$
is a fractional ideal of $\O_S$.
Finally, given that $\mathfrak{b}$ and $\O_S^\perp$ are fractional ideals, it is easy to show that
$\mathfrak{b}^\perp=\mathfrak{b}^{-1}\O_S^\perp$.
\qed

\vspace{1ex}


\subsection{The action of the units is ergodic}

\begin{definition*}
Let $G$ be a compact topological group with normalized Haar measure $\mu$.
An automorphism $\sigma:G\to G$ is ergodic if
$\sigma(E)=E$ implies
$\mu(E)=0$ or $\mu(E)=1$
for every measurable set $E$.
\end{definition*}

\begin{lemma}[Halmos]\label{L:halmos}
A (continuous) automorphism of a
compact abelian group $G$ is ergodic if and only if the induced automorphism
on the character group $G^\vee$ has no finite orbits (other than the trivial one).
\end{lemma}

\noindent\textbf{Proof.}
The proof is a one page argument using Pontryagin duality and Fourier series.
See~\cite{halmos} for the details.
\qed

\begin{lemma}\label{L:meta.ergodic}
If $u\in\U_S$ is not a root of unity, then the automorphism of $\T$
given by $[\xi]\mapsto u[\xi]$ is ergodic.
\end{lemma}

\noindent\textbf{Proof.}
We will use Lemma~\ref{L:halmos}.
Since any character of $\T$ may be viewed as a character of $\overline{K}_S$
that is trivial on $\mathfrak{a}$, we will consider characters of $\overline{K}_S$.
One checks that the action of multiplication by $\U_S$
induces the action $u\phi_\xi=\phi_{u\xi}$
on $\overline{K}^\vee_S$.

Now let $n$ denote a nonzero integer.
Using duality, we have
$
  u^n\phi_\xi=\phi_\xi
  \Rightarrow
  \phi_{u^n\xi}=\phi_\xi
  \Rightarrow
  u^n\xi=\xi
  \Rightarrow
  (u^n-1)\xi=0
  \Rightarrow
  u^n=1 \text{ or } \xi=0
$.
By hypothesis, $u$ is not a root of unity.
Hence the only solution to
$u^n\chi=\chi$ is when $\chi$ is the trivial character.
\qed

%


\subsection{Convergence of subgroups}

\begin{definition*}
Let $G$ be a locally compact abelian group.  For a subgroup $H$ of $G$, we define
the annihilator of $H$ in $G^\vee$ to be
$$
A(G^\vee, H)=\{\chi\in G^\vee\mid \chi(H)=1\}
\,.
$$
\end{definition*}

\begin{lemma}[Berend]\label{L:berend}
Let $G$ be a compact abelian metric group.  A sequence $G_n$ of closed subgroups of $G$
satisfies $G_n\to G$ in the Hausdorff metric if and only if for every nonzero
$\chi\in G^\vee$ we have $\chi\notin A(G^\vee, G_n)$ for sufficiently large $n$.
\end{lemma}

\noindent\textbf{Proof.}
The proof is half a page and uses the Haar measure and integral
on the groups involved.
See~\cite{berend:multi} for the details.
\qed

\begin{lemma}\label{L:chartheory}
Suppose $L$ is a subgroup of $\overline{K}_S$.  Let $u\in\U_S$.
If $u^n\overline{\pi(L)}\not\to\T$, then there exists nonzero
$\alpha\in K$ and an increasing sequence $n_k\in\Z^+$ such that
$\phi(u^{n_k}\alpha L)=1$ for all
$k$.
\footnote{Actually $\alpha$ lies in the fractional ideal $\mathfrak{a}^\perp$,
but we won't need this.}

\end{lemma}

\noindent\textbf{Proof.}
Suppose
$u^n\overline{\pi(L)}\not\to\T$.
Then Lemma~\ref{L:berend} says that there exists a nonzero $\chi\in\T^\vee$ 
and an increasing sequence 
$n_k\in\Z^+$ such that \mbox{$\chi(u^{n_k}\overline{\pi(L)})=1$} for all $k$.
Viewing $\chi$ as a character on $\overline{K}_S$
and using duality (Lemma~\ref{L:dual}), we know there exists a nonzero $\alpha\in\mathfrak{a}^\perp$
so that $\chi([\xi])=\phi_\alpha(\xi)=\phi(\alpha\xi)$ for all $\xi\in\overline{K}_S$.
This leads to $\phi(u^{n_k}\alpha L)=1$ for all $k$.
\qed

\section{Proof of Proposition~\ref{P:3}}\label{S:Proof2}


\begin{lemma}\label{L:ergodic}
Let $U$ be a finite index subgroup of $\U_S$,
and $\Lambda$ be a closed $U$-invariant subset of $\T$ with nonempty interior.
If $\#S\geq 2$, then 
$\Lambda=\T$.
\end{lemma}

\noindent\textbf{Proof.}
First observe that $\Lambda$ has nonzero measure because it has nonempty interior.
The group $U$ is of finite index in $\U_S$ and therefore
$\rank(U)=\rank(\U_S)=\#S-1\geq 1$, which implies that $U$ contains a unit
which is not a root of unity.
Now Lemma~\ref{L:meta.ergodic} implies that $\Lambda$ is dense in $\T$,
giving the result.~\qed


\vspace{1ex}

In order to prove the proposition, we first give a construction and a lemma.
For the remainder of this section, let $M$ be a $\U_S$-minimal subset of $\T$.
It suffices to show that $M-M=\T$ implies $M=\T$,
as clearly $\T$ is not $\U_S$-minimal.
Hence we assume that $M-M=\T$.
We will write $\xi,\eta$
for an element of $\overline{K_S}$ as well as the corresponding element of $\T$;
that is, we will drop the brackets from the expressions $[\xi]$,$[\eta]$.

\vspace{1ex}

\noindent\textbf{Construction.}
Define $U^{(n)}=(\U_S)^{n!}$ so that $U^{(n)}\subseteq (\U_S)^n$
and
$
  \U_S=U^{(1)}\supseteq U^{(2)}\supseteq U^{(3)}\supseteq\dots\,;
$
choose a sequence of subsets
$
M=M^{(1)}\supseteq M^{(2)}\supseteq\dots
$
so that $M^{(k)}$ is $U^{(k)}$-minimal.
Finally, define the set $M^{\infty}=\cap_k M^{(k)}$;
observe that $M^\infty$
is closed and nonempty since $\T$ is compact.

\begin{lemma}
Given $\xi\in K$, we have $\xi+\eta\in M$ for all $\eta\in M^\infty$. 
\end{lemma}

\noindent\textbf{Proof.}
Let $\xi\in K$.
First we show that there exists $\eta'\in M^\infty$ so that $\xi+\eta'\in M$.
Since $\Orb(\xi)$ is finite, there exists $N\in\Z^+$ such that
$(\U_S)^N\xi=\{\xi\}$ and hence $U^{(N)}\xi=\{\xi\}$.
For ease of notation, set $U'=U^{(N)}$ and $M'=M^{(N)}$.
Since $U'$ has finite index in $\U_S$, we have
$\U_S/ U'=\{a_1U',a_2U',\dots,a_\ell U'\}$
for $a_k\in\U_S$ with $a_1=1$.

We define the closed sets $\Lambda_i=M-a_iM'$
for $i=1,\dots,\ell$.
We observe
$\cup_{i=1}^\ell a_i M' =M$
as the former set is closed and $\U_S$-invariant
and clearly contained in the latter set.
Since $M-M=\T$ by hypothesis, this leads to
$\cup_{i=1}^\ell \Lambda_i=\T$.
It is now easy to see that $\Lambda_j$ must have nonempty interior
for some $j$.
Since $\Lambda_j$ is closed and $U'$-invariant, 
Lemma~\ref{L:ergodic} gives $\Lambda_j=\T$.

It follows that there exists $\eta\in a_j M'$ such that
$a_j\xi+\eta\in M$ and therefore
$\xi+ a_j^{-1} \eta=a_j^{-1}(a_j\xi+\eta)\in M$.
It is plain that $\eta':=a_j^{-1}\eta\in M'$.
We have shown that there exists $\eta'\in M'$ such that $\xi+\eta'\in M$.
Now observe that $U'(\xi+\eta')\subseteq M$, and
since $\overline{U'\eta'}=M'$ and $U'\xi=\{\xi\}$,
we have $\xi+\eta\in M$ for all $\eta\in M'$.
\qed

\vspace{1ex}

\noindent\textbf{Proof of Proposition~\ref{P:3}.}
Fix $\eta\in M^\infty$.  We will show that $M-\eta=\T$ from which
$M=\T$ immediately follows.  Since $K$ is dense in $\T$,
it suffices to show that $K\subseteq M-\eta$.
Let $\xi\in K$ be arbitrary.
The previous lemma says that $\xi+\eta\in M$.
The result follows.
\qed


\section{Proof of Proposition~\ref{P:2}.}\label{S:Proof1}

%

The following standard result in algebraic number theory
will be helpful.
If one was forced to attach names to it, the following might be called the
Dirichlet--Minkowski--Hasse--Chevalley Unit Theorem.

\namedthm{$S$-Unit Theorem}
{
For every $w\in S$ there exists $\eps\in\U_S$ such that $|\eps|_v<1$ for all $v\in S$ with $v\neq w$.
Moreover, choosing $\eps_w$ as above for each $w\in S$ yields a set $\{\eps_w\}_{w\in S}$
which, after any one element is discarded, forms
an independent set of units (modulo torsion) and generates a finite index subgroup of $\U_S$;
in particular $\rank(\U_S)=\#S-1$.
}

The following lemma allows us to locate points that
``live in a single component''.

\begin{lemma}\label{L:force}
Suppose $\#S\geq 2$.
Let $N$ be a closed $\U_S$-invariant subset of $\overline{K}_S$
that contains $0$ as a non-isolated point.
For each $w\in S$, the set $N\cap K_w$ contains a nonzero point.
\end{lemma}

\noindent\textbf{Proof.}
By hypothesis, there is a sequence $\xi_n\in N$, $\xi_n\neq 0$, with $\xi_n\to 0$.
We will write $\xi_n=(\xi_{n,v})_{v\in S}$.
By the $S$-Unit Theorem there exists a unit $u\in\U_S$ such that
$|u|_v<1$ for all $v\neq w$, and hence $C:=|u|_w>1$.
Define
$$
  \mathcal{A}
  =
  \left\{(\alpha_v)_{v\in S}\in\overline{K}_{S}
  :
  |\alpha_w|_w\geq 1
  \,,\;
  |\alpha_v|_v\leq C \;\;\forall v\in S  
  \right\}
  \,.
$$
For all sufficiently large $m$
we have $|\xi_{m,v}|_v\leq 1$ for all $v\in S$ and hence
there exists a $j_m\in\Z^+$ such that
$u^{j_m}\xi_m\in\mathcal{A}$.
Since $\mathcal{A}$ is compact there is a limit point $\eta$
of this sequence; $\eta\in\mathcal{A}$ and hence $\eta\neq 0$.
Since $N$ is $\U_S$-invariant and closed we have $\eta\in N$.
Finally, for all $v\neq w$ we have $|u|_v<1$ which implies
$|u^{j_m}\xi_{m,v}|_v\to 0$ and hence $\eta\in K_w$.
\qed

\subsection{$K$ has a real embedding}

Given what we have shown up to this point, it is now quite easy to establish
Proposition~\ref{P:2} in the case where $K$ has a real embedding.
This makes use of the following fact:

\begin{lemma}\label{L:dense}
Suppose $K\subseteq\R$ is a number field.
If $\rank(\U_S)\geq 2$, then $\U_S$ is dense in $\R$.
\end{lemma}

This result is well-known and not hard to establish, but we prove it here for
the sake of completeness and
also because it motivates what we do in the general case.  We will need the following well-known
result in Diophantine approximation (see, for example, \cite{cassels:book}).

\namedthm{Kronecker's Theorem}
{
Let $\alpha_1,\dots,\alpha_n\in\R$.
Then
$$
  \{(m\alpha_1,\dots,m\alpha_n)\mid m\in\Z\}
$$
is dense in $\R^n/\Z^n$
iff $1,\alpha_1,\dots,\alpha_n$ are linearly independent over $\Q$.
}

Lemma~\ref{L:dense} follows immediately from:

\begin{lemma}
Let $a,b\in\R^+$.  Suppose $a$ and $b$ are multiplicatively independent.
Then $\{a^n b^m\mid n,m\in\Z\}$ is dense in $\R^+$.
\end{lemma}

\noindent\textbf{Proof.}
Taking the logarithm to the base $a$ of $a^n b^m$ gives
$m+n\alpha$ where $\alpha=\log b/\log a$.
Thus it suffices to show that $\{m+n\alpha\mid m,n\in\Z\}$ is dense in $\R$.
But since $a,b$ are multiplicatively independent, we know that
$\alpha$ is irrational.  Thus $\{n\alpha\mid n\in\Z\}$ is dense in $\R/\Z$
by Kronecker's Theorem.
The result follows.
\qed

\vspace{2ex}

\noindent\textbf{Proof of Proposition~\ref{P:2} when $K$ has a real embedding.}\\
Set $\tilde{N}=\pi^{-1}(N)$.
Then $\tilde{N}$ is a closed $\U_S$-invariant subset of $\overline{K}_S$
that contains $0$ as a non-isolated point.
Let $w$ be a real place and apply Lemma~\ref{L:force} to $\tilde{N}$.
This gives an element $x\in\R=K_w\subseteq\overline{K}_S$
such that $x\in\tilde{N}$, $x\neq 0$.
Lemma~\ref{L:dense} tells us that $\U_S$ is dense in $\R$
and hence $\{ux\mid u\in\U_S\}$ is dense in $\R$;
it follows that $\tilde{N}$ contains $\R$.
Now Lemma~\ref{L:prelim.help} gives $\tilde{N}=\overline{K}_S$
and hence $N=\T$.
\qed

\vspace{1ex}

At this junction, we point out that we have completely justified
Theorem~\ref{T:3}, and hence all the results of \S\ref{S:intro},\ref{S:main.idea},\ref{S:euclidean},
in the case where $K$ has a real embedding.
In particular, this is enough to establish Corollary~\ref{C:2}.
However, there is more work to be done to establish our results
in the case where $K$ is totally complex.
We don't seem to get any additional
mileage out of the assumption that all the embeddings are complex,
so we will simply work with number fields
that have at least one complex embedding.


\subsection{$K$ has a complex embedding}

We recall the following standard definition.
\begin{definition*}
We call a number field $K$ a CM-field if either of the two
equivalent conditions are satisfied:
\begin{enumerate}
\item
$K$ is a totally complex quadratic extension of a totally real field.
\item
There is a subfield $F$ of $K$ with $\rank(\U_F)=\rank(\U_K)$.
\end{enumerate}
\end{definition*}
The equivalence of the two definitions follows from Dirichlet's Unit Theorem
(or the $S$-Unit Theorem with $S=S_\infty$).  We write $K^+$ for the
maximal totally real subfield of $K$.  In the case where $K$ is a CM-field we have
$[K:K^+]=2$.

\begin{lemma}\label{L:CMkey}
Suppose $K\subseteq\C$ is a number field with $K\not\subseteq\R$.  If $K$ is not a CM-field,
then there exists $u\in\U$ such that $u^n\notin\R$ for all nonzero $n\in\Z$.
\end{lemma}

\noindent\textbf{Proof.}
Suppose that for every $u\in\U$ there exists $n\in\Z^+$ such that $u^n\in\R$.
It follows that there must exists $N\in\Z^+$ such that
$\U^N\subseteq\R$.  If $K\not\subseteq\R$, this implies $\Q(\U^N)\neq K$ which 
forces $K$ to be a CM-field.
\qed

\begin{lemma}\label{L:CMkey2}
Suppose $K\subseteq\C$ is a number field with $K\not\subseteq\R$.  If $K$ is a CM-field
and $S$ contains a finite prime that splits in $K/K^+$,
then there exists $u\in\U_S$ such that $u^n\notin\R$ for all nonzero $n\in\Z$.
\end{lemma}

\noindent\textbf{Proof.}
Let $\mathfrak{P}$ be a finite prime in $S$ that splits in $K/K^+$.
Let $h$ denote the class number of $K$.  Then define $u\in\O$
by $(u)=\mathfrak{P}^h$.  It is plain that $u\in\U_S$ since $v(u)=0$
for all $v\notin S$.  By way of contradiction, suppose
$u^n\in\R$ for some nonzero $n\in\Z$.
Then we would have $u^n\in K^+$ and $(u^n)=\mathfrak{P}^{hn}$ in $K$.
Since $\mathfrak{P}$ lies above two distinct primes in $K^+$,
this is impossible.
\qed


\vspace{1ex}

In what follows, we will write $[x]$ to denote the floor of $x$, and
write $x=[x]+\{x\}$ so that $\{x\}$ denotes the
fractional part of $x$.
We will also use the notation $\|x\|=\inf_{y\in\Z}|x-y|$.

\begin{lemma}
Suppose $\alpha,\beta\in\R$ and $\alpha\notin\Q$.  Then there exists $r,s\in\Z$ with
$r>0$ such that
$\{(m\alpha,m\beta)\mid m\in\Z\}$ is dense in
$\{(rt,st)\mid t\in\R\}$
when they are both viewed as subsets of $\R^2/\Z^2$.
\end{lemma}

\noindent\textbf{Proof.}
If $\{1,\alpha,\beta\}$ is linearly independent over $\Q$ then the result follows from
Kronecker's Theorem (with $n=2$).  Otherwise we have $a\alpha+b\beta+c=0$
with $a,b,c\in\Z$, not all zero; we must have 
$b\neq 0$ lest we contradict the fact that $\alpha\notin\Q$,
and, without loss of generality, we may assume that $b>0$.
Pick $t\in\R$
and let $\eps>0$ be given.
Pick $\delta>0$ so that $|a|\delta,|b|\delta<\eps$.
Applying Kronecker's Theorem (with $n=1$) we may
choose $m\in\Z$ so that $\|m\alpha-t\|<\delta$.
It follows that $\|mb\alpha-bt\|<|b|\delta<\eps$.
Also we have $b\beta=-a\alpha-c$ which implies
$
  mb\beta=-a(m\alpha-t)-at-mc
  \,.
$
Therefore $\|mb\beta+at\|<|a|\delta<\eps$.
It follows that $(mb\alpha,mb\beta)$ is $\eps$-close to $(bt,-at)$ in $\R^2/\Z^2$.
\qed

\vspace{1ex}

The following result says that in our situation the closure of
$\U_S$ contains a nice spiral or concentric circles.
It is a little complicated to state, but it plays the same role as
Lemma~\ref{L:dense}.

\begin{lemma}\label{L:spiral}
Suppose $K\subseteq\C$ and $K\not\subseteq\R$.
Suppose $K$ is not a CM-field or
$S$ contains a finite prime that splits in $K/K^+$.
If $\#S\geq 3$ then
either:
\begin{enumerate}
\item
$\overline{\U}_S\supseteq\{z^t\mid t\in\R\}$
where
$z\in\C\setminus\R$, $|z|>1$
\item
$\overline{\U}_S\supseteq\{x^n z^t \mid t\in\R,\, n\in\Z\}$
where
$z\in\C\setminus\R$, $|z|=1$,
$x\in\R$, $x>1$
\end{enumerate}
\end{lemma}

\noindent\textbf{Proof.}
First, suppose there exists $u\in\U_S$ with $|u|=1$ which is not a root of unity.
In this case, $\{u^m\mid m\in\Z\}$ is dense in the unit circle.
Using the $S$-Unit Theorem we may choose $v\in\U_S$ with $|v|>1$.
We see that $\{u^m v^n\mid m,n\in\Z\}$ is dense in
$\{|v|^n u^t \mid t\in\R,\,n\in\Z\}$.
In this case, conclusion 2 holds with $z=u$ and $x=|v|$.
Hence we may assume that no elements of $\U_S$ other than
roots of unity are unimodular.

Since $\rank(\U_S)=\#S-1\geq 2$, we
know that $\U_S$ contains two independent units $u$ and $v$.
Write $u=|u|e^{2\pi i\theta}$ and $v=|v|e^{2\pi i\varphi}$, where
$|u|,|v|\neq 1$.
Given our hypotheses, we may assume that
$\theta\notin\Q$
(see Lemmas~\ref{L:CMkey} and~\ref{L:CMkey2}).
Without loss of generality, we may assume $|u|>1$
by replacing $u$ with $u^{-1}$ if necessary.

Set $\alpha=\log |v|/\log |u|$ and
$\beta=\varphi-\alpha\theta$.
Observe that $\alpha$ is irrational; if it was the case that
$\alpha=a/b$, then $u^av^{-b}$ would be a unimodular
unit which is not a root of unity.
Choose $r,s$ as in the previous lemma.
Set $z=|u|^{r}e^{2\pi i(r\theta+s)}$.
%
Choose $\delta>0$ small enough so that
$\{rt\}=r t$ for all $t\in[0,\delta]$.
Since $\overline{\U}_S$ is a multiplicative group, to prove the lemma
it suffices to show that
$\overline{\U}_S\supseteq\{z^t\mid 0< t< \delta\}$.

Fix $t\in(0,\delta)$.
We construct sequences $n_k$ and $m_k$
so that $u^{n_k} v^{m_k}\to z^t$.
By our choice of $r,s$, there is a sequence $m_k$ so that
$(m_k\alpha, m_k\beta)$ converges to $(rt, st)$ in $\R^2/\Z^2$;
it follows that
$\{m_k\alpha\}\to\{r t\}$ and
that $m_k\beta$ converges to $st$ in $\R/\Z$.
Set $n_k=-[m_k\alpha]$ so that $n_k+m_k\alpha\to r t$.
It follows that $|u|^{n_k}|v|^{m_k}\to |u|^{r t}$.
Now observe that
$$
  n_k\theta+m_k\varphi
  =
  \{m_k\alpha\}\theta + m_k\beta
  \,,
$$
which converges (modulo $1$) to $r t\theta+s t$.
It follows that 
$$
u^{n_k}v^{m_k}\to  |u|^{r t}e^{2\pi i(r\theta+s)t}=z^t
\,.
\;\;\;\text{\qed}
$$

\vspace{2ex}

\noindent\textbf{Proof of Proposition~\ref{P:2}.}
Set $\tilde{N}=\pi^{-1}(N)$.
Then $\tilde{N}$ is a closed $\U_S$-invariant subset of $\overline{K}_S$
that contains $0$ as a non-isolated point.
Pick a complex place of $K$ and apply Lemma~\ref{L:force} to $\tilde{N}$.
(Since we have already proved the result when $K$ has a real place,
we may certainly assume that $K$ has a complex place.)
This gives a nonzero element $a\in\C\subseteq\overline{K}_S$
such that $a\in\tilde{N}$.
In what follows, distances between sets and convergence of sets
will always be measured using the standard Hausdorff distance.

\vspace{2ex}

\noindent\textbf{Claim 1:}
There exists a sequence of (compact) arcs $A_n$
and line segments $L_n$
which lie in $\C\subseteq\overline{K}_S$
with the following properties:
$$
A_n\subseteq\tilde{N}
\,,\;
d(L_n,A_n)\to 0
\,,\;
\length(L_n)\to\infty
$$

We apply Lemma~\ref{L:spiral} and obtain one
of two possible conclusions (see the statement of the lemma).
First, we assume that conclusion 1 holds.
(When conclusion 2 holds, the proof will be similar.)
We have that $\overline{\U}_S$
contains $\{z^t\mid t\in\R\}$ for some $z\in\C\setminus\R$ with $|z|>1$.
Therefore $\tilde{N}$ contains the spiral $\{a z^t\mid t\in\R\}$ in $\C$.
The arcs $A_n$ we construct will be subarcs of this spiral
and are therefore all automatically contained in $\tilde{N}$.
Let $\delta_n>0$ be a sequence of real numbers
with $\delta_n\to 0$ to be chosen later.
Define the arc:
$$
  A_n=\{a z^t
  :
  t\in[n,n+\delta_n]\}
  \,.
$$
Let $L_n$ denote the corresponding line segment which
has the same endpoints.
Namely,
$$
  L_n=
  \left\{
  a z^n\left[1+\lambda(z^{\delta_n}-1)\right]
  :
  \lambda\in[0,1]
  \right\}
  \,.
$$
For $\lambda\in[0,1]$ we write $t=n+\delta_n\lambda$
and, using calculus, we obtain:\footnote{$\Log$ will denote the principal branch of the logarithm.}
\begin{eqnarray*}
  A_n(\lambda)-L_n(\lambda)
  &=&
  az^n\left[
  (z^{\lambda\delta_n}-1)
  -
  \lambda(z^{\delta_n}-1)
  \right]
  \\
  &=&
  az^n
  \left[
    \frac{\Log^2 z}{2}
    \lambda(\lambda-1)\delta_n^2+O(\delta_n^3)
  \right]
\end{eqnarray*}
Thus there are constants $C_1,C_2>0$ such that
for $n$ sufficiently large, we have:
\begin{eqnarray*}
  d(A_n,L_n)
  &\leq&
  C_1 |z|^n\delta_n^2\\[1ex]
  \length(L_n)
  &=&
  |a| |z|^n|z^{\delta_n}-1|
  \\
  &\geq&
  C_2 |z|^n\delta_n
\end{eqnarray*}
Choose $\delta_n>0$ with
$|z|^n\delta_n^2\to 0$ but $|z|^n\delta_n\to\infty$
so that
$d(A_n,L_n)\to 0$ and $\length(L_n)\to\infty$.
This completes the proof of the claim in this case.

Now we suppose that conclusion 2 of Lemma~\ref{L:spiral} holds,
so that $\tilde{N}$ contains $\{ax^n z^t\mid t\in\R,\,n\in\Z\}$;
here $z\in\C\setminus\R$, $|z|=1$, $x\in\R$, $x>1$.
This time we use:
$$
  A_n=\{a x^n z^t : t\in[0,\delta_n]\}
$$
$$
  L_n=
  \left\{
  a x^n \left[1+\lambda(z^{\delta_n}-1)\right]
  :
  \lambda\in[0,1]
  \right\}
$$
For $\lambda\in[0,1]$ we write $t=\delta_n\lambda$
and find:
$$
  A_n(\lambda)-L_n(\lambda)
  =
  ax^n\left[
  (z^{\lambda\delta_n}-1)
  -
  \lambda(z^{\delta_n}-1)
  \right]
$$
Hence there are constants $C_1,C_2>0$ such that
for $n$ sufficiently large, we have:
$$
  d(A_n,L_n)
  \leq
  C_1 x^n\delta_n^2
  \,,
  \quad
  \length(L_n)
  \geq
  C_2 x^n\delta_n
$$
As before we choose $\delta_n$ appropriately and the 
claim follows.

\vspace{2ex}

\noindent\textbf{Claim 2:}
There is a line $L\subseteq\C\subseteq\overline{K}_S$
passing through the origin
such that $[\xi]+\pi(L)\subseteq N$
for some $\xi\in\overline{K}_S$.

\vspace{2ex}

We can think of each line segment $L_n$ (given by the previous claim) as a triple
$(x_n,y_n,z_n)\in\C\times S^1\times \R^+$
which represents the midpoint, direction, and length
of the segment.  Namely,
$$
  L_n=\{x_n+t y_n:-z_n\leq 2t\leq z_n\}
  \,.
$$
(The choice of $y_n\in S^1$ in this representation
is not unique, but this won't affect the argument.)
By passing to a subsequence, we may assume
that $\pi(x_n)\to [\xi]$ for some $[\xi]\in\T$ and $y_n\to y$ for some $y\in S^1$;
we have $z_n\to\infty$ by what we have already shown.
Therefore, for each $t\in\R$ we have
$\pi(x_n+t y_n)\to [\xi]+\pi(ty)$.

Let $L$ denote the line corresponding to the triple $(0,y,\infty)$
which passes through the origin in the direction of $y$;
namely,
$$
  L=\{t\,y\mid t\in\R\}
  \,.
$$
We show that $[\xi]+\pi(L)\subseteq N$.
Let $\eta\in[\xi]+\pi(L)$ be arbitrary and $\eps>0$ be given.
For $n$ sufficiently large, we have
$$
  d(\eta,N)\leq d(\eta,\pi(A_n))\leq d(\eta,\pi(L_n))+d(\pi(L_n),\pi(A_n))<\eps
  \,.
$$
Since $N$ is closed, we obtain $\eta\in N$.
This proves the claim.

\vspace{2ex}

\noindent\textbf{Claim 3:}
There is a unit $u\in\U_S$ such that $u^n\overline{\pi(L)}\to\T$.

\vspace{2ex}

Pick $u\in\U_S$ so that $u^n\notin\R$ for all $n\in\Z^+$
(see Lemmas~\ref{L:CMkey} and~\ref{L:CMkey2}).
By way of contradiction, suppose  $u^n\overline{\pi(L)}\not\to\T$.
In light of Lemma~\ref{L:chartheory}
there exists a nonzero $\alpha\in K$ and a strictly increasing sequence of positive integers $n_k$
so that
$\phi(u^{n_k}\alpha L)=1$ for all $k$.
Fix an arbitrary $k\in\Z^+$.
We have
$\phi(u^{n_k}\alpha t y)=1$ for all $t\in\R$.
Since $y\in\C\subseteq\overline{K}_S$ and the local character is
$\phi_w(z)=e^{2\pi i (z+\overline{z})}$ (see Equation~\ref{E:phi}),
this leads to $2\Re(u^{n_k}t y')\in\Z$ for all $t\in\R$
where we define $y':=\rho_w \alpha y\in\C\subseteq\overline{K}_S$.
Because $y'\neq 0$, it follows that $\Re(u^{n_k}y')=0$.
Now we see that $u^{n_2-n_1}\in\R$.
This contradiction proves the claim.

\vspace{2ex}

By Claim~2 and the fact that $N$ is closed, we have
$[\xi]+\overline{\pi(L)}\subseteq N$.
By Claim~3, we have $u^n \overline{\pi(L)}\to\T$.
Choose a subsequence so that
$u^{n_k}[\xi]\to[\eta]$ for some $[\eta]\in\T$.
Therefore
$u^{n_k}([\xi]+\overline{\pi(L)})=u^{n_k}[\xi]+u^{n_k}\overline{\pi(L)}\to[\eta]+\T$.
Since $N$ is closed and $\U_S$-invariant,
we conclude $[\eta]+\T\subseteq N$ which implies $N=\T$.
\qed


\section{Proof of Proposition~\ref{P:4}}\label{S:CM}

Suppose Theorem~\ref{T:3} holds except in the case where $K$ is a CM-field and
distinct primes in $S$ lie over distinct primes in $K^+$.  We will show it holds
in the remaining cases.

Let $K$ be a CM-field with totally real subfield $K^+$.
We assume that distinct primes in $S$ lie over distinct primes in $K^+$.
Let $S^+$ denote the set of all primes in $K^+$
(including the infinite ones) lying under primes in $S$.
Given our hypotheses, no finite primes of $S$ are split
in $K/K^+$ and $\#S=\#S^+$.
For ease of notation, we will write $\U_S=\U_{K,S}$ and
$\U_{S^+}=\U_{K^+,S^+}$. Now we proceed to the proof proper.

\vspace{2ex}

\noindent\textbf{Proof of Proposition~\ref{P:4}.}
Since $\U_{S^+}$ is contained in $\U_{S}$,
it suffices to prove that every  nonempty closed $\U_{S^+}$-invariant
subset of $\T$ contains torsion elements.

First we consider what happens locally.
Choose $v\in S^+$.  By our hypothesis, there is exactly one prime $w\in S$
lying above $v$.
In this case,
we have the inclusion $K^+_v\subseteq K_w$
(with $[K_w:K^+_v]=2$)
and the isomorphism
$K_w/K^+_v\simeq K^+_v$.
(Let $\{1,\theta\}$ be a basis for $K/K^+$; then $\{1,\theta\}$
is also a basis for $K_w/K^+_v$ and the
aforementioned isomorphism is the one that sends the
coset represented by $x+y\theta$ to the element $y$.)
The multiplication action of $\U_{S^+}$ on $K_w$ induces an action
on $K^+_v$ via mapping $K_w\rightarrow K_w/K^+_v\simeq K^+_v$;
one checks that this is just the usual multiplication action so that in
subsequent arguments we won't be dealing with two different actions.
%

In light of the inclusions of local fields discussed above, we have
an inclusion
$\overline{K}^+_{S^+}\subseteq\overline{K}_S$ and an isomorphism
$\overline{K}_S/\overline{K}^+_{S^+}\simeq \overline{K}^+_{S^+}$.
If we define
$\mathfrak{a}^+=\mathfrak{a}\cap K^+$ and
$\T^+=\overline{K}^+_{S^+}/\mathfrak{a}^+$,
then this leads to an exact sequence
of compact abelian groups:
\begin{equation}\label{E:exact}
  0\to\T^+\to\T\to\T/\T^+\to 0
\end{equation}
We have that $\U_{S^+}$ acts on $\T$ and $\T^+$ and this leads 
to an action of $\U_{S^+}$ on $\T/\T^+$.
Moreover, Theorem~\ref{T:3} applies to $K^+$ and hence
every nonempty $\U_{S^+}$-invariant subset of $\T^+$ contains
torsion elements.
We will show that the same holds for $\T/\T^+$ and hence for $\T$.

The situation is summarized by the following commutative diagram:
  $$
  \xymatrix
  {
  & 0 \ar@{->}[d] & 0 \ar@{->}[d] & 0 \ar@{->}[d] \\
  0 \ar@{->}[r] & \mathfrak{a}^+ \ar@{->}[r]  \ar@{->}[d] &\mathfrak{a} \ar@{->}[r]  \ar@{->}[d] & \mathfrak{a}/\mathfrak{a}^+ \ar@{->}[r] \ar@{->}[d] & 0\\
  0 \ar@{->}[r] & \overline{K}^+_{S^+} \ar@{->}[r]  \ar@{->}[d] &\overline{K}_S \ar@{->}[r]  \ar@{->}[d] & \overline{K}_S / \overline{K}^+_{S^+} \ar@{->}[r]  \ar@{->}[d] & 0\\
  0 \ar@{->}[r] & \T^+ \ar@{->}[r]  \ar@{->}[d] &\T \ar@{->}[r]  \ar@{->}[d] &  \T/\T^+ \ar@{->}[r] \ar@{->}[d] & 0\\
  & 0 & 0 & 0 \\
  }
  $$
The isomorphism
$\overline{K}_S/\overline{K}^+_{S^+}\simeq\overline{K}^+_{S^+}$
carries $\mathfrak{a}/\mathfrak{a}^+$ to a fractional
ideal $\mathfrak{b}$ and we obtain the isomorphic exact sequences:
$$
\xymatrix
{
   0
   \ar@{->}[r]
   &
   \mathfrak{a}/\mathfrak{a}^+ \ar@{->}[r]
   \ar@{->}[d]\ar@<0.7ex>@{}[d]|-*@{~}   
   &
   \overline{K}_S / \overline{K}^+_{S^+}
   \ar@{->}[r]
   \ar@{->}[d]\ar@<0.7ex>@{}[d]|-*@{~}   
   &
   \T/\T^+
   \ar@{->}[r]
   \ar@{->}[d]\ar@<0.7ex>@{}[d]|-*@{~}   
   &
   0\\
   0 \ar@{->}[r]  & \mathfrak{b} \ar@{->}[r] & \overline{K}^+_{S^+} \ar@{->}[r] & \T_*^+ \ar@{->}[r] & 0 \\
}
$$
In light of previous comments, the induced action of $\U_{S^+}$
on $\T^+_*= \overline{K}^+_{S^+}/\mathfrak{b}$ is
just the usual multiplication action.
Here we have written $\T^+_*$ to emphasize that it is potentially different
from $\T^+$ since the quotient is by a different ideal.  In any case,
Theorem~\ref{T:3} applied to $K^+$ tells us that
every nonempty closed $\U_{S^+}$-invariant subset of $\T^+_*$
contains torsion elements; hence the same holds for $\T/\T^+$.

Finally, we use what we have shown to complete the proof.
Reconsidering the exact sequence~(\ref{E:exact}),
we note that desired result holds for both $\T^+$ and $\T/\T^+$.
Let $N$ be a nonempty closed $\U_{S^+}$-invariant subset of $\T$.
Applying the result for $\T/\T^+$,
we have that $\pi(N)\subseteq\T/\T^+$ contains torsion elements
which implies\footnote{Here $\pi$ denotes the projection $\pi:\T\to\T/\T^+$
which is different than the previous usage of this notation.}
that $mN\cap \T^+\neq\emptyset$ for some $m\in\Z^+$.
Since the set $mN\cap \T^+$ has the property in question,
we apply the result for $\T^+$ to conclude that $mN\cap\T^+$ contains
torsion elements.  It follows that $N$ contains torsion elements.
\qed


\section*{Acknowledgements}
This author would like to thank Hendrik Lenstra for his
helpful suggestions and encouragement,
and Mary Flahive and Tom Schmidt for many
helpful discussions.


\bibliographystyle{plain}
\bibliography{myrefs}

\vspace{2ex}

{\footnotesize
\noindent
Kevin J. McGown\\
Oregon State University\\
Corvallis, Oregon, USA\\[2ex]
\emph{Current Address:}  Ursinus College, 601 E. Main St., Collegeville, Pennsylvania, USA\\
\emph{Email:}  kmcgown@ursinus.edu
}

\end{document}